\documentclass{amsart}

\usepackage{epsfig}
\usepackage{amsthm}
\usepackage{amssymb}
\usepackage{amsmath}
\usepackage{amscd}

\newtheorem {thm}   {Theorem}

\numberwithin{thm}{section}

\newtheorem {lem}[thm]    {Lemma}
\newtheorem {prop}[thm]{Proposition}
\theoremstyle{definition}
\newtheorem{de}[thm]{Definition}
\theoremstyle{remark}
\newtheorem{rem}[thm]{Remark}

\newtheorem {cor}[thm]{Corollary}

\newtheorem{Notation}[thm]{Notation}
\expandafter\chardef\csname pre amssym.def
at\endcsname=\the\catcode`\@ \catcode`\@=11
\def\undefine#1{\let#1\undefined}
\def\newsymbol#1#2#3#4#5{\let\next@\relax
 \ifnum#2=\@ne\let\next@\msafam@\else
 \ifnum#2=\tw@\let\next@\msbfam@\fi\fi
 \mathchardef#1="#3\next@#4#5}
\def\mathhexbox@#1#2#3{\relax
 \ifmmode\mathpalette{}{\m@th\mathchar"#1#2#3}%
 \else\leavevmode\hbox{$\m@th\mathchar"#1#2#3$}\fi}
\def\hexnumber@#1{\ifcase#1 0\or 1\or 2\or 3\or 4\or 5\or 6\or 7\or 8\or
 9\or A\or B\or C\or D\or E\or F\fi}

\font\teneufm=eufm10 \font\seveneufm=eufm7 \font\fiveeufm=eufm5
\newfam\eufmfam
\textfont\eufmfam=\teneufm \scriptfont\eufmfam=\seveneufm
\scriptscriptfont\eufmfam=\fiveeufm

\catcode`\@=\csname pre amssym.def at\endcsname

\def \Z {{\mathbb Z}}
\def \C {{\mathbb C}}

\def \Q {{\mathbb Q}}

\def \wt{{\rm wt}}

\def \Res{{\rm Res}}

\def \End{{\rm End}}

\def \mod{{\rm mod}}
\def \<{\langle}
\def \>{\rangle}

\def \a{\alpha }

\def \g{\gamma}
\def \b{\beta }

\def \qed{\mbox{ $\square$}}
\def \pf {\noindent {\bf Proof:} \,}
\def \cg{\chi_g}
\def \cg'{\chi'_g}

\def \o{\otimes}
\def \d{\delta}

\def \r{\rho}

\def \l{\lambda }

\begin{document}






\title[$C_2$-cofiniteness of the vertex operator algebra $V_L^+$]{$C_2$-cofiniteness of the vertex operator algebra $V_L^+$ when $L$ is a rank one lattice}
\author[Gaywalee Yamskulna]{Gaywalee Yamskulna}
\address{Department of Mathematics, SUNY, Binghamton, NY 13902 and}
\address{Mathematical Science Research Institute, Berkeley, CA 94720}
\email{gail@math.binghamton.edu}
\date{\today}
\maketitle
\begin{abstract}
Let $L$ be a rank one positive definite even lattice. We prove
that the vertex operator algebra (VOA) $V_L^+$ satisfies the $C_2$
condition. Here, $V_L^+$ is the fixed point sub-VOA of the VOA
$V_L$ associated with the automorphism lifted from the -1 isometry
of $L$.
\end{abstract}
\bigskip
\section{Introduction}

A vertex operator algebra (VOA) $V$ is said to be {\em $C_2$
cofinite }if the subspace $\{ u_{-2}v|u,v\in V\}$ has finite
codimension in $V$. This is often call the $C_2$ condition and was
first introduced in [Z] by Zhu who used it to prove the convergence of the
trace function of a certain kind of VOA-modules. This seemingly
abstract condition is satisfied by most known VOAs, and has played
an important role in the representation theory of VOAs and the
study of the structure of their modules. In particular, under this
condition, it is possible to establish the existence of twisted
modules (see [DLM1]). Furthermore, for holomorphic VOAs,
 it was shown in [DLM1] that the $C_2$ condition implies the
uniqueness of twisted $V$-modules. It was also used, in [KL], to
show that every irreducible admissible twisted $V$-module is an
ordinary twisted $V$-module. Most recently, the $C_2$ condition
was shown to imply the finiteness of the generating sets of a
given VOA and its modules (see [GN, Bu]).

A VOA $V$ is called rational if any $V$-module is completely
reducible. Understanding the representation theory of rational
VOAs is one of the major problems in the field and it seems that
here too the $C_2$ condition will play an important role. In fact,
it has been conjectured that rationality and $C_2$-cofiniteness
are equivalent. Indeed, in [L], Li showed that any regular vertex
operator algebra satisfies the $C_2$ condition.

If the central charge is less than 1, the representations of
rational VOAs have been completely understood because the sub-VOAs
generated by the Virasoro elements have been completely classified
(see [DMZ, W]). Hence, the first nontrivial case, is when the
central charge is 1. Let $L$ be a rank one positive definite even
lattice. It is well known that the corresponding VOA $V_L$ has an
order 2 automorphism $\theta$ which is induced from the -1
isometry of the lattice. The $\theta$-invariant sub-VOA $V_L^+$ is
a simple VOA (see [DM]). It has been conjectured that every
rational VOAs of central charge 1 is of the form $V_L$, $V_L^+$
and $V^G_{L_2}$, where $L_2$ is a root lattice of type $A_1$, $G$
is a finite subgroup of $SO(3)$ of type $E$, and $V^G_{L_2}$ is a
$G$-invariant sub-VOA of $V_{L_2}$. The representation theory of
$V_L$ is completely understood. In fact, $V_L$ is rational (see
[Bo, FLM, D1, DLM3]). Therefore, in order to characterize the
rationality of vertex operator algebras with central charge 1, one
has to understand $V_L^+$ and $V^G_{L_2}$. To this end, Dong and
Nagatomo classified the irreducible modules for $V_L^+$ (see
[DN2]). In this paper, we take a step towards achieving this goal
by {\em showing that $V_L^+$ satisfies the $C_2$ condition.}

Let $L=\Z \a$ be an even lattice with a non-degenerate integral
bilinear form $\<\cdot, \cdot\>$ such that $\<\a,\a\>=2k$, where
$k$ is a positive integer. When $k=1,2$, $C_2$-cofiniteness is a
consequence of the fact that $V_L^+$ is isomorphic to VOAs which
are known to be $C_2$ cofinite. More specifically, when $k=1$,
$V_L^+$ is isomorphic to a particular lattice VOA (see [DG]), and
when $k=2$ it is isomorphic to $L(1/2,0)\o L(1/2,0)$ (see [DGH]).
The proof that these VOAs are $C_2$ cofinite is available in
[DLM1].

Here we consider the case {\em $k\geq 3$} and establish $C_2$
cofiniteness in the following way. First we show that
$V_L^++C_2(V_L^+)$ is generated by $M(1)^++V_L^+(1)+C_2(V_L^+)$
where $V_L^+(1)=M(1)^+\o (e^{\a}+e^{-\a})\oplus M(1)^- \o
(e^{\a}-e^{-\a})$. Then we use information about the bases of
$M(1)^+$ and $V_L^+(1)$ to show that $(M(1)^++C_2(V_L^+))/
C_2(V_L^+)$ and $(V_L^+(1)+C_2(V_L^+))/C_2(V_L^+)$ have finite
dimension. This implies that $V_L^+$ is $C_2$ cofinite.

The paper is organized as follows. In section $2$, we review the
definitions of a vertex operator algebra, its automorphisms, and
its twisted modules. We discuss the definition the cofiniteness $C_n$
condition and the algebra of $V/C_2(V)$. We also recall the
construction of $A_n(V)$ for a nonnegative integer $n$. In section
$3$, we review the construction of vertex operator algebras
$V_L^+$, its irreducible modules and discuss the cofiniteness
$C_2$ condition of $V_L^+$ when $k=1,2$. In section $4$, we show
that $V_L^+/C_2(V_L^+)$ is spanned by the subspace
$M(1)^++V_L^+(1)+C_2(V_L^+)$. Finally, we show that $V_L^+$ is
$C_2$ cofinite in section $5$.

\bigskip

\section{\bf{Vertex operator algebras and $C_2$ condition}}

We review definitions of a vertex operator algebra, its
automorphisms, and its twisted modules. We also discuss the definition
of the cofiniteness $C_n$ condition and the algebra of $V/C_2(V)$.
We recall from [DLM2] the construction of $A_n(V)$ for a
nonnegative integer $n$.

\smallskip

For a vector space $W$, we let $W[[z, z^{-1}]]$ be the space of
$W$-valued formal series in arbitrary integral powers of $z$.

\smallskip

\begin{de}\label{d2.1} {[Bo, FLM, FHL]} A {\it vertex operator algebra} (or VOA) is a
${\Z}$-graded vector space
$$V=\oplus_{n\in{\Z}}V_n;$$ such that
\begin{eqnarray}
\dim\,V_n &<&\infty\ \ \mbox{ and}\\
V_n&=& 0\ \ \mbox{ if}\ \ n\ \ \mbox{is\ \ sufficiently\ \ small.}
\end{eqnarray}
Moreover, there is a linear map
\begin{equation}
\begin{array}{l}
V \to (\mbox{End}\,V)[[z,z^{-1}]]\\
v\mapsto Y(v,z)=\displaystyle{\sum_{n\in{\Z}}v_nz^{-n-1}}\ \ \ \
(v_n\in\mbox{End}\,V)
\end{array}
\end{equation}
and with two distinguished vectors ${\bf 1}\in V_0$, $\omega \in
V_2$ satisfying the following conditions for $u, v \in V$:
\begin{eqnarray}
u_nv &=&0\ \ \ \ \ \mbox{for}\ \  n\ \ \mbox{sufficiently
large};\label{e2.2}\\
Y({\bf 1},z)&=&1;\label{e2.3}\\
Y(v,z){\bf 1}\in V[[z]]&\mbox{and}&\lim_{z\to 0}Y(v,z){\bf 1}=v;
\end{eqnarray}
\begin{equation}\label{jac}
\begin{array}{c}
\displaystyle{z^{-1}_0\delta\left(\frac{z_1-z_2}{z_0}\right)
Y(u,z_1)Y(v,z_2)-z^{-1}_0\delta\left(\frac{z_2-z_1}{-z_0}\right)
Y(v,z_2)Y(u,z_1)}\\
\displaystyle{=z_2^{-1}\delta \left(\frac{z_1-z_0}{z_2}\right)
Y(Y(u,z_0)v,z_2)}
\end{array}
\end{equation}
(Jacobi identity) where $\delta(z)=\displaystyle{\sum_{n\in
{\Z}}z^n}$ is the algebraic formulation of the $\delta$-function
at 1, and all binomial expressions are to be expanded in
nonnegative integral powers of the second variable;
\begin{equation}\label{e2.6}
[L(m),L(n)]=(m-n)L(m+n)+\frac{1}{12}(m^3-m)\delta_{m+n,0}(\mbox{rank}\,V)
\end{equation}
for $m, n\in {\Z},$ where
\begin{equation}
L(n)=\omega_{n+1}\ \ \ \mbox{for}\ \ \ n\in{\Z}, \ \ \
\mbox{i.e.},\ \ \ Y(\omega,z)=\sum_{n\in{\Z}}L(n)z^{-n-2}
\end{equation}
and
\begin{eqnarray}
& &\mbox{rank}\,V\in {\Q};\\
& &L(0)v=nv=(\mbox{wt}\,v)v \ \ \ \mbox{for}\ \ \ v\in V_n\
(n\in{\Z}); \label{3.40}\\
& &\frac{d}{dz}Y(v,z)=Y(L(-1)v,z). \label{3.41}
\end{eqnarray}
\end{de}
We denote the vertex operator algebra just defined by
$(V,Y,\bf{1},\omega)$ (or briefly, by $V$). The series $Y(v,z)$
are called {\it vertex operators.}

\smallskip

\begin{de}
A VOA $V$ is of {\em CFT} type if $V=\bigoplus_{n\geq 0}V_n$ and
$V_0= \C {\bf 1}$.
\end{de}

\smallskip

\begin{de}
An {\it automorphism} of  $(V,Y,\bf{1},\omega)$ is a linear map
$g$: $V\to V$  satisfying
\begin{eqnarray*}
gY(v,z)g^{-1}&=&Y(gv,z),\ \ v\in V,\\
g\bf{1}&=&\bf{1},\\
g{\omega}&=&\omega.
\end{eqnarray*}
\end{de}

\smallskip

For $g$, an  automorphism of the VOA $V$ of order $T$, we denote
the decomposition of $V$ into eigenspaces with respect to the
action of $g$ as $V=\bigoplus_{r=0}^{T-1}V^r$ where $V^r=\{v\in
V|gv=e^{2\pi ir/T}v\}$. For a vector space $W$, we denote the
space of $W$-valued formal series in arbitrary complex powers of
$z$ by $W\{z\}$.
\begin{de}
A {\em weak $g$-{\it twisted} $V$-module} $M$ is a vector space
equipped with a linear map
\begin{equation}
\begin{array}{l}
V\to (\mbox{End}\,M)\{z\}\label{map}\\
v\mapsto\displaystyle{ Y_M(v,z)=\sum_{n\in\Q}v_nz^{-n-1}\ \ \
(v_n\in\End M)}
\end{array}
\end{equation}
satisfying  axioms analogous to (\ref{e2.2}), (\ref{e2.3}) and
(\ref{jac}). To describe these, we let $u\in V^r$, $v\in V$ and
$w\in M$. Then
\begin{eqnarray}
& &Y_M(u,z)=\sum_{n\in r/T+\Z}u_nz^{-n-1};
\label{1/2}\\
& &u_nw=0\ \ \ \mbox{for}\ \ \ n\ \ \ \mbox{sufficiently\
large};\label{ds1}\\
& &Y_M({\bf 1},z)=1;
\end{eqnarray}
\begin{equation}\label{jacm}
\begin{array}{c}
\displaystyle{z^{-1}_0\delta\left(\frac{z_1-z_2}{z_0}\right)
Y_M(u,z_1)Y_M(v,z_2)-z^{-1}_0\delta\left(\frac{z_2-z_1}{-z_0}\right)
Y_M(v,z_2)Y_M(u,z_1)}\\
\displaystyle{=z_2^{-1}\left(\frac{z_1-z_0}{z_2}\right)^{-r/T}\delta\left(\frac{z_1-z_0}{z_2}\right)
Y_M(Y(u,z_0)v,z_2)}.
\end{array}
\end{equation}
 \end{de}
We denote this module by $(M,Y_M),$ or briefly by $M$. Equation
(\ref{jacm}) is called the {\em twisted Jacobi identity}. If $g$
is the identity element, this reduces to the definition of a weak
$V$-module and (\ref{jacm}) is the untwisted Jacobi identity .

\smallskip

\begin{de}\label{d1} An ordinary {\em $g$-twisted $V$-module} is a
weak $g$-twisted $V$-module $M$ which carries a $\C$-grading
induced by the spectrum of $L(0)$. Then $$M=\bigoplus_{\lambda\in
\C}M_{\lambda}$$ where $M_{\lambda}=\{w\in M|L(0)w=\lambda w\},$
dim $M_{\lambda}<\infty.$ Moreover, for fixed $\lambda, M_{{n\over
T}+{\lambda}}=0 $ for all small enough integers $n.$
\end{de}

\smallskip

\begin{de}\label{d2} An {\em admissible $g$-twisted $V$-module} is a
weak $g$-twisted $V$-module $M$ which carries a ${1\over T}\Z_{+}$
grading $M=\bigoplus_{n\in {1\over T}\Z_{+}}M(n)$ satisfying the
following condition:$$v_{m}M(n)\subset M(n+\wt v-m-1)$$ for
homogeneous $v\in V.$ Here, $\Z_+$ is the set of nonnegative
integers.
\end{de}

\smallskip

\begin{rem} The notion of admissible $g$-twisted $V$-module here is equivalent
to the notion of a module in [Z] when $g$ is the identity element.
\end{rem}

\smallskip

\begin{lem}{[DLM4]}
Any $g$-twisted $V$-module is an admissible $g$-twisted
$V$-module.
\end{lem}
So, there is a natural identification of the category of
$g$-twisted $V$-modules with a sub-category of the category of
admissible $g$-twisted $V$-modules.

\smallskip

\begin{de} $V$ is called rational if every
admissible $V$-module is a direct sum of irreducible admissible
$V$-modules.
\end{de}

\smallskip

\begin{de} $V$ is called regular if every weak $V$-module is a direct sum of irreducible ordinary $V$-modules.
\end{de}

\begin{cor} If $V$ is a regular VOA then $V$ is a rational VOA.
\end{cor}

\smallskip

We now recall the definition of the $C_n$-condition and discuss
about the algebra of $V/C_2(V)$.

\begin{de}
For a VOA $V$, we define $$C_n(V)=\{ \ \ v_{-n}u\ \ | \ \ v, u\in
V\ \ \}.$$ $V$ is said to satisfy {\em the cofiniteness $C_n$
condition} if $V/C_n(V)$ is finite dimensional.
\end{de}

\smallskip

\begin{rem}
For the case $n=2$, it was introduced by Zhu (see [Z]).
\end{rem}

\smallskip

The following lemma is a consequence of a definition \ref{d2.1}.
\begin{lem}\label{l2.5}\ \

1) $L(-1)u\in C_2(V)$ for all $u\in V$.

2) $v_{-n}u\in C_2(V)$ for all $u, v\in V$, and $n\geq 2$.
\end{lem}
\pf Part 1), it follows from the fact that $L(-1)u =u_{-2}{\bf 1}$
for all $u\in V$.

Part 2), it follows from the fact that $(L(-1)^m
v)_{-2}=(m+1)!v_{-m-2}$.\qed

\smallskip

\begin{rem} $C_2(V)=\sum_{n\geq {2}}C_n(V).$
\end{rem}

\smallskip

Next, we discuss the algebra of $V/C_2(V)$. For a VOA $V$, we
consider the $(-1)^{st}$ product
\begin{eqnarray*}
V\times V &\rightarrow & V\\
a\times b & \mapsto& a\cdot b =a_{-1}b.\\
\end{eqnarray*}
By using the Jacobi Identity, we obtain the following.

\smallskip

\begin{thm}{[Z]} \ \

1) $C_2(V)$ is an ideal of $V$ with respect to the $(-1)^{st}$
product.

2) $V/C_2(V)$ is a commutative associative algebra under
$(-1)^{st}$ product.
\end{thm}

\smallskip

\begin{thm}{[GN]}
Let $U=\{ u^i\}_{i\in I}$ be a set of homogeneous elements in $V$
which are representatives of a basis of $V/C_2(V)$. Then $V$ is
spanned by elements of the form
$$u^{i_1}_{-n_1}u^{i_2}_{-n_2}...u^{i_l}_{-n_l}{\bf 1}$$ where
$n_1>n_2>...>n_l>0$ and $u^{i_j}\in U$.
\end{thm}

\smallskip

\begin{thm}{[GN]}
Let $V$ be a VOA of {\em CFT}-type. If $V/C_2(V)$ is finite
dimensional, then $V/C_n(V)$ is finite dimensional for $n\geq 2$.
\end{thm}

\smallskip

We review the associative algebra $A_n(V)$ constructed in [DLM2].
\begin{de} Let $n$ be a nonnegative integer. We define $$O_n(V)=\<\ \ u \circ _n v,\ \ L(-1)u+L(0)u\ \ |\ \ \mbox{homogeneous}\ \ u, v\in V\>.$$ Here $u\circ _n v=\Res _z Y(u,z)v{(1+z)^{ \wt u+n}\over z^{2n+2}}.$
Also, we define a product $*_n$ on $V$ for $u, v$ as above;
$$u*_nv=\sum_{m=0}^n(-1)^m\left(\begin{array}{c}
                m+n\\n
                \end{array}\right)\Res _z Y(u,z)v{(1+z)^{ \wt u +n}\over z^{n+m+1}}.$$
\end{de}
Set $A_n(V)=V/O_n(V)$.
\begin{thm}{[DLM2]}\ \

1) $O_n(V)$ is a 2-sided ideal of $V$ under $*_n$.

2) $A_n(V)$ is an associative algebra under $*_n$ with the
identity ${\bf 1}+O_n(V)$. Moreover, $\omega + O_n(V)$ is a
central element of $A_n(V)$.

\end{thm}

\smallskip

\begin{rem} $A_0(V)$ is an $A(V)$ in [Z] and it was first introduced by Zhu (see [Z]).
\end{rem}

\smallskip

\begin{thm}{[Bu]}
If $V$ is a simple VOA that satisfies the $C_2$ condition, then
the associative algebra $A_n(V)$ is finite dimensional for all
$n\geq 0$.
\end{thm}

\bigskip

\section{\bf{Vertex operator algebra $V_L^+$}}

In this section, we briefly review the construction of $V_L^+$,
its irreducible modules and discuss the $C_2$ condition of the VOA
$V_L^+$ when $k=1,2$.

We are working in the setting of [FLM, DL1]. Let $L=\Z\a$ be an
even lattice with a non-degenerate integral bilinear form
$\<\cdot,\cdot \>$ such that $\<\a,\a\>=2k$. Here, $k$ is a
positive integer. Set ${\mathbf{h}}=\C\otimes_{\Z}L$ and extend
the form  $\< \cdot,\cdot \>$ from $L$ to $\mathbf{h}$ by
$\C$-bilinearity. Let $\widehat{\bf h}={\bf h}\o \C [t,
t^{-1}]\oplus \C c$ be the affinization of ${\bf h}$. Therefore,
$\widehat{\bf h}$ is a Lie algebra with commutator relations:

\begin{eqnarray*}
[h\o t^m, h'\o t^n]&=& m\d_{m+n,0}\< h, h'\> c \ \ \text{for}\ \ h,h'\in {\bf h}, m,n\in \Z;\\
{[ } \widehat{\bf h}, c {]} &=& 0.\\
\end{eqnarray*}
For $h\in {\bf h}, n\in \Z$, we use the notation $h(n)$ to denote
$h\o t^n$. Set $${\widehat{\bf h}}^+= {\bf h}\o t\C[t],\ \
{\widehat{\bf h}}^-={\bf h}\o t^{-1}\C[t^{-1}].$$ Hence,
$\widehat{\bf h}^+$, $\widehat{\bf h}^-$ are abelian subalgebra of
$\widehat{\bf h}$.

For a Lie algebra ${\bf g}$, we let $U({\bf g})$ be the universal
enveloping algebra of $\bf g$. Consider the induced $\widehat{\bf
h}$-module
$$M(1)=U(\widehat{\bf h})\o _{U({\bf h}\o \C[t]\oplus \C c)} \C,$$
where ${\bf h}\o \C[t]$ acts trivially on $\C$ and $c$ acts on
$\C$ as multiplication by 1. We set $\C[L]$ to be a group algebra
of $L$ with a basis $\{ e^{\b}|\b\in L\}$. Let $z$ be a formal
variable and $h\in {\bf h}$. We define actions of $h$ and $z^{h}$
on $\C[L]$ in the following ways:
\begin{eqnarray*}
h\cdot e^{\b}&=&\< h, \b\>e^{\b};\\
z^h\cdot e^{\b}&=&z^{\< h,\b \>}e^{\b}.
\end{eqnarray*}

We let $$V_L=M(1)\o \C[L].$$ $L, \widehat{\bf h}, z^h (h\in {\bf
h})$ act naturally on $V_L$ by acting on either $M(1)$ or $\C[L]$
as indicated above. Next, we shall define the vertex operator
$Y(v,z)$ for $v\in V_{L}$. For $\b\in L, \g,\a _1,..., \a _k\in
{\bf{h}}, n_1,..., n_k\in \Z(n_i> 0)$, we set
\begin{eqnarray*}
\g(z)&=&\sum_{n\in \Z}\g (n)z^{-n-1},\ \ \text{where}\ \ \g(n)=\g\o t^n,\\
Y(e^{\b},z)&=&exp\left(\sum_{n=1}^{\infty}\frac{\b(-n)}{n}z^n\right)exp\left(- \sum_{n=1}^{\infty}\frac{\b(n)}{n}z^{-n}\right)e^{\b} z^{\b},\\
Y(v,z)&=&:\left({1\over{(n_1-1)!}}\left({d\over{dz}}\right)^{n_1-1}\a_1(z)\right)...\left({1\over{(n_k-1)!}}\left({d\over{dz}}\right)^{n_k-1}\a_k(z)\right)Y(e^{\b},z):,
\end{eqnarray*}
where $v=\a _1(-n_1)...\a_k(-n_k)\o e^{\b}$.
 We use a normal ordering procedure, indicated by open colons, which
signify that the enclosed expression is to be reordered if
necessary so that all the operators $\g(n),(\g\in {\bf{h}},
n<0),e^{\b} \in {L}$ are to be placed to the left of all operators
$\a(n), z^{\a},(\a\in {\bf{h}}, n\geq 0)$ before the expression is
evaluated. This gives a well defined linear map
$$\begin{array}{lll}
V_L&\to &(\End\,V_{L })[[z, z^{-1}]]\\
v&\mapsto&\displaystyle{ Y(v,z)=\sum_{n\in\Z}v_nz^{-n-1}\ \ \
(v_n\in \End\,V_L)}.
\end{array}$$

\smallskip

\begin{thm}\label{thex1} The space $(V_L, Y,{\bf 1}, \omega)$ is a simple VOA
with $\omega={1\over 4k}\a(-1)^2$ and ${\bf 1}=1\o 1$ ( see [Bo,
FLM]).
\end{thm}

\smallskip

Let $L^{\circ}=\{ x\in {\bf h}|\< x, L\>\subset \Z\}$ be the dual
lattice of $L$. Hence, $L^{\circ}=\frac{1}{2k}L$ and
$L^{\circ}=\bigcup_{i=0}^{2k-1}(L+\frac{i}{2k}\a)$ is a coset
decomposition of $L^{\circ}$ with respect to $L$.

\smallskip

\begin{thm}\ \

1) The irreducible $V_L$-modules are
$$V_{L+\frac{i}{2k}\a}=M(1)\o \C[L+\frac{i}{2k}\a]$$ where $i=0,..., 2k-1$ (see
[D1]).

2) $V_L$ is a rational vertex operator algebra (see [DLM3]).

\end{thm}

We set $$V_{L^{\circ}}= \bigoplus_{i=0}^{2k-1}V_{L+{i\over
2k}\a}.$$ Let $\theta$ be the linear automorphism of $V_{L^\circ}$
such that $\theta (u\o e^{\g})=\theta (u)\o e^{-\g}$ for $u\in
M(1)$ and $\g\in L^{\circ}$. The action of $\theta$ on $M(1)$ is
given by
$$\theta (\a _1(n_1)...\a _k(n_k))=(-1)^k \a _1(n_1)...\a
_k(n_k).$$ The restriction of $\theta$ to $V_L$ is a VOA
automorphism. Let $M$ be an $\theta$-stable subspace of $V_{L}$.
We denote ${\pm}1$ eigenspaces by $M^{\pm}$, respectively. Note
that $M(1)^+$ is a vertex operator subalgebra of $V_L^+$.

\begin{thm}\label{pex5}\ \

1) ${V_L}^{+}$ is a simple VOA (see [DM]).

2) Set $J={1\over 4k^2}\a(-1)^4{\bf 1}-{1\over k}\a(-3)\a(-1){\bf
1}+{3\over 4k}\a(-2)^2{\bf 1}$. Then $V_L^+$ is generated by
$\omega$, $e^{\a}+e^{-\a}$, and $J$ (see [DG]).

\end{thm}

\smallskip

Let $T_1, T_2$ be $\Z _2$-irreducible modules such that $\a +2L$
acts as scalars 1 and -1, respectively. We set
$$\widehat{\bf{h}}[-1]={\bf h}\o t^{1/2}\C[t,t ^{-1}]\oplus \C c.
$$
Hence, $\widehat{\bf h}[-1]$ is a Lie algebra with the commutator
relation given by

\begin{eqnarray*}
[h\o t^m, h'\o t^n]&=& m\d_{m+n,0}\< h, h'\> c \ \ \text{for}\ \ h,h'\in {\bf h}, m,n\in \Z+\frac{1}{2};\\
{[ } \widehat{\bf h}, c {]} &=& 0.\\
\end{eqnarray*}
Moreover, $\widehat{\bf h}[-1]^+={\bf h}\o t^{1/2}\C[t]\oplus \C
c$ is a subalgebra of $\widehat{\bf h}[-1]$. We define the action
of $\widehat{\bf h}[-1]^+$ on $T_i$ $(i=1,2)$ by the actions
$$h\o t^{1/2+n}\cdot u =0 \ \ \text{and}\ \ c\cdot u=u\ \ ( h\in
{\bf h},n\in \Z_{\geq 0}, u\in T_i).$$

\smallskip

For $i=1,2$, we set $$V_L^{T_i}=U(\widehat{{\bf h}}[-1])\o
_{U(\widehat{{\bf h}}[-1]^+)}T_i.$$
\smallskip

\begin{lem} \ \

1) $V_L^{T_i}$ is the irreducible $\theta$-twisted $V$-module (see
[FLM], [DL2]).

2) $V_L^{T_i}( i=1,2)$ are all irreducible $\theta$-twisted
$V_L$-modules (see [D2]).

\end{lem}

\smallskip

We define a linear operator $\theta$ on $V_L^{T_i}$: for $\a_i\in
{\bf h}, n_i\in \frac{1}{2}+\Z, t\in T_i$
$$\theta(\a_1(-n_1)...\a_s(-n_s)\o
t=(-1)^s\a_1(-n_1)...\a_s(-n_s)\o t. $$ In fact, $\theta$ is an
automorphism of $V_L^{T_i}$. Then we have the decomposition
$$V_L^{T_i}=(V_L^{T_i})^+\oplus (V_L^{T_i})^-.$$

\smallskip

\begin{lem}{[FLM]} $(V_L^{T_i})^{\pm}$ are irreducible
$V_L^+$-modules for $i=1,2$.
\end{lem}

\smallskip

\begin{thm}{[DN2]}
$\{V_L^{\pm}, V_{L+\frac{\a}{2}}^{\pm},V_{L+\frac{r\a}{2k}}
(V_L^{T_i})^{\pm}|i=1,2, 1\leq r\leq k-1\}$ is the set of all
inequivalent irreducible $V_L^+$-module.
\end{thm}

\smallskip

\begin{thm}{[DN2]} $A(V_L^+)$ is a semisimple algebra.
\end{thm}

\smallskip

\begin{lem} When $k=1,2$, $V_L^+$ satisfies the $C_2$-condition.
\end{lem}
\pf For $k=1$, $V_L^+$ is isomorphic to the lattice VOA $V_{L'}$
where $L'$ is a rank one positive definite lattice spanned by $\b$
whose square length is 8 (see. [DG]). Since $V_{L'}$ satisfies the
$C_2$-condition (see [DLM1]), this implies that $V_L^+$ satisfies
the $C_2$-condition when $k=1$.

For $k=2$, $V_L^+$ is isomorphic to $L(1/2,0)\o L(1/2,0)$ where
$L(1/2,h)$ is the irreducible highest weight module for the
Virasoro algebra with central charge 1/2 and highest weight $h$
(see [DGH]). It was proved in [DLM1] that $L(1/2,0)\o L(1/2,0)$
satisfies the $C_2$ condition. Therefore, $V_L^+$ satisfies the
$C_2$ condition. \qed

\bigskip

\section{\bf{ A spanning set of $V_L^+/C_2(V_L^+)$}}

For the rest of this paper, we {\em assume $k\geq 3$}. In this
section, we show that $V_L^+/C_2(V_L^+)$ is generated by the
subspace $M(1)^++V_L^+(1)+C_2(V_L^+)$. Here, $$V_L^+(1)=M(1)^+\o
(e^{\a}+e^{-\a})\oplus M(1)^- \o (e^{\a}-e^{-\a}). $$

\smallskip

Denote by $\Z_{\geq 0}$ the set of nonnegative integers, by
$\Z_{>0}$ the set of positive integers. We recall that

\[{1\over n!}\left({d\over dz}\right)^n\a(z)=\sum_{j\geq 0}\left(\begin{array}{c}
                -j-1\\n
                \end{array}\right)\a(j)z^{-j-n-1}+\sum_{j\leq -n-1}\left(\begin{array}{c}
                -j-1\\n
                \end{array}\right)\a(j)z^{-j-n-1}\]\\
where $\a(z)=\sum_{m\in\Z}\a(m)z^{-m-1}$. For $m\in \Z_{>0}$, we
set
\begin{eqnarray*}
E^m&=&e^{m\a}+e^{-m\a}\ \ \mbox{, and}\\
F^m&=&e^{m\a}-e^{-m\a}.\\
\end{eqnarray*}
For convenience, we also set
\begin{eqnarray*}
E&=&E^1 \ \ \mbox{, and}\\
F&=&F^1.\\
\end{eqnarray*}

\begin{rem} $\a(0)E^m=2km F^m$, and $\a (0) F^m=2km E^m$.
\end{rem}

\smallskip

We set
\begin{eqnarray*}
\exp\left(\sum_{n=1}^{\infty}{x_n\over n}z^n\right)&=& \sum_{j=0}^{\infty}p_j(x_1,x_2,...)z^j\\
&=&\sum_{j=0}^{\infty}p_j(x)z^j.\\
\end{eqnarray*}
Then $p_j(x)$ are Schur polynomials. For $m\in\Z_{>0}$, we let $$V_L^+(m)= M(1)^+\o E^m +M(1)^- \o
F^m.$$ We will show that if $m$ is even, then $V_L^+(m)\subset
M(1)^++C_2(V_L^+)$, and if $m$ is odd, then $V_L^+(m)\subset
V_L^+(1)+C_2(V_L^+)$.

\begin{lem}\label{sl1}
For $m\in \Z_{>0}$, $E^{2m}\in M(1)^+ +C_2(V_L^+)$.
\end{lem}
\pf Since
\begin{eqnarray*}
Y(E^m, z)E^m&=& Y(e^{m\a},z)e^{m\a}+ Y(e^{-m\a},z)e^{m\a}+Y(e^{m\a},z)e^{-m\a}+Y(e^{-m\a},z)e^{-m\a}\\
&=&\sum_{j=0}^{\infty}p_j(m\a)\o e^{2m\a}z^{2m^2k+j}+\sum_{j=0}^{\infty}p_j(-m\a)z^{-2m^2k+j}\\
& & +\sum_{j=0}^{\infty}p_j(m\a)z^{-2m^2k+j}+\sum_{j=0}^{\infty}p_j(-m\a)\o e^{-2m\a}z^{2m^2k+j}\\
&=& \sum_{j=0}^{\infty}(p_j(m\a)\o e^{2m\a}+p_j(-m\a)\o e^{-2m\a})z^{2m^2k+j}\\
& & + \sum_{j=0}^{\infty}(p_j(m\a)+p_j(-m\a))z^{-2m^2k+j},\\
\end{eqnarray*}
we have $\left(E^m\right)_{-2m^2k-1}E^m=
E^{2m}+p_{4m^2k}(-m\a)+p_{4m^2k}(m\a)$. Hence $E^{2m}\in  M(1)^+
+C_2(V_L^+). \ \ \qed$

\medskip

\begin{lem}\label{sl2}
For $m\in \Z_{\geq 0}$, $E^{2m+1}\in V_L^+(1) +C_2(V_L^+)$.
\end{lem}
\pf Assume $m\neq 0$. Since
\begin{eqnarray*}
Y(E^m, z)E^{m+1}&=&\sum_{j=0}^{\infty}(p_j(m\a)\o e^{(2m+1)\a}+p_j(-m\a)\o e^{-(2m+1)\a})z^{2m(m+1)k+j}\\
& & + \sum_{j=0}^{\infty}(p_j(-m\a)\o e^{\a}+p_j(m\a)\o e^{-\a})z^{-2m(m+1)k+j},\\
\end{eqnarray*}
it follows that
\begin{equation*}
\left(E^m\right)_{-2km(m+1)-1}E^{m+1}= E^{2m+1}+
p_{4km(m+1)}(-m\a)\o e^{\a}+p_{4km(m+1)}(m\a)\o e^{-\a}.
\end{equation*}
Therefore, we have $E^{2m+1}\in V_L^+(1) +C_2(V_L^+)$. \qed

\smallskip

\begin{lem}\label{sl3} {[DN2]}
For $m,n \geq 1$, $$(\a(-n)\a(-1){\bf 1})_{-1}E^m=
2km(n+(-1)^{n-1})\a(-n-1)F^m+\a(-n)\a(-1)E^m.$$
\end{lem}

\smallskip

\begin{cor}\label{sc1}\ \

1) If $m$ is even, then we have
\begin{equation}\label{e14}
2km(n+(-1)^{n-1})\a(-n-1)F^m+\a(-n)\a(-1)E^m\in  M(1)^+
+C_2(V_L^+).
\end{equation}

2) If $m$ is odd, then we have
\begin{equation}
2km(n+(-1)^{n-1})\a(-n-1)F^m+\a(-n)\a(-1)E^m\in V_L^+(1)
+C_2(V_L^+).
\end{equation}
\end{cor}
\pf These follow immediately from Lemmas \ref{sl1}, \ref{sl2},
\ref{sl3} and the fact that $C_2(V_L^+)$ is an ideal of $V$ under
$(-1)^{st}$ product. \qed

\medskip

Let $u=\a(-n_1)...\a(-n_r)v$, ($n_i>0,$ $v=E^m$ or $F^m$).
Following [DN2], we say that {\em an element $u$ has length $r$
with respect to $\a$ and we write $l_{\a}(u)=r$}. In general, if
$u$ is a linear combination of such vectors $u^i$ we define the
length of $u$ to be the maximal length among $l_{\a}(u^i)$.

\smallskip

\begin{lem} {[DN2]}\label{sl4.5} Let $m\in \Z_{>0}$.

1) Let $n_1,n_2,..., n_r\in \Z_{>0}$ with $r$ even. Then
\begin{equation*}
(\a(-n_1)\a(-n_2)...\a(-n_r){\bf 1})_{-1}
E^m=\a(-n_1)\a(-n_2)...\a(-n_r)E^m +u
\end{equation*}
where $l_{\a}(u)<r$, $u\in V_L^+(m)$.

\ \

2) Let $n_1,..., n_r\in \Z_{>0}$ with $r$ odd. Then
\begin{equation*}
(\a(-n_1)\a(-n_2)...\a(-n_{r-1}){\bf 1})_{-1}(\a(-n_r)
F^m)=\a(-n_1)\a(-n_2)...\a(-n_r)F^m +u
\end{equation*}
where $l_{\a}(u)<r$, $u\in V_L^+(m)$.
\end{lem}

\smallskip

\begin{lem} Let $m, n\in \Z_{>0}$. We have
\begin{equation}\label{se1}
n\a(-n-1) F^m\equiv -m\a(-n)\a(-1)E^m \ \ \mod\ \ C_2(V_L^+).
\end{equation}
\end{lem}
\pf This Lemma follows from the Lemma \ref{l2.5}, and the fact
that
\begin{eqnarray*}
L(-1)\a(-n)F^m &=& {1\over 4k}\{ 2\a(-n-1)\a(n)+ 2\a(-1)\a(0)\}\a(-n)F^m\\
 & =& n\a(-n-1)F^m+m\a(-n)\a(-1)E^m. \ \ \qed\\
\end{eqnarray*}

\smallskip

\begin{lem}\label{l4.7}
Let $m, n\in \Z_{>0}$. If $m$ is even, then we have
$$\a(-n-1)F^m\equiv 0 \ \ \mod \ \ M(1)^+ +C_2(V_L^+).$$
\end{lem}

\pf By combining equations (\ref{e14}), and (\ref{se1}) together,
we obtain that:
\begin{equation*}
2km^2(n+(-1)^{n-1})\a(-n-1)F^m\equiv n\a(-n-1)F^m\ \ \mod \ \
M(1)^+ + C_2(V_L^+).
\end{equation*}
Consequently, we have
\begin{equation*}
(2km^2(n+(-1)^{n-1})-n)\a(-n-1)F^m\in   M(1)^+ + C_2(V_L^+).
\end{equation*}

Assume that $(2km^2(n+(-1)^{n-1})-n)=0$.

\underline{case 1}: $n$ is an even integer. Then we have
$n={2km^2\over 2km^2 -1}$. Recall that for any integer $p\neq 0$,
if $p|p+1$ then $p=1$ or $p=-1$. Therefore, either $2km^2=2$ or
$2km^2=0$. This contradicts with the fact that $2km^2\geq 6$.

\underline{case 2:} $n$ is an odd integer. Then $n=-{2km^2\over
2km^2-1}$. This is impossible because $n$ is a positive integer.

Thus, $(2km^2(n+(-1)^{n-1})-n)\neq 0$, and $\a(-n-1)F^m\in M(1)^+
+ C_2(V_L^+)$. \qed

\smallskip

\begin{lem}
If $m$ is odd, then  $\a(-n-1)F^m\in V_L^+(1) + C_2(V_L^+)$ for
all $n\geq 1$.
\end{lem}

\smallskip

\begin{lem}\label{l4.9}
For $m\in \Z_{>0}$, $\a(-1)F^m\in C_2(V_L^+)$.
\end{lem}

\pf This follows from the fact that $L(-1)E^m =m\a(-1)F^m$. \qed

We recall that for $m\in \Z_{>0}$, $V_L^+(m)=M(1)^+\o E^m+
M(1)^-\o F^m$.

\begin{lem}\label{l4.12}
For any even positive integer $m$, $V_L^+(m)\subset
M(1)^++C_2(V_L^+).$
\end{lem}

\pf Let $u=\a(-n_1)...\a(-n_r)(e^{m\a}+(-1)^re^{-m\a})\in
V_L^+(m)$. We will prove this lemma by using an induction on $r$.
When $r=0$, it follows immediately from Lemma \ref{sl1}. When
$r=1$, it follows from Lemmas \ref{l4.7}, \ref{l4.9}. Set
\begin{equation}
v=
\begin{cases}
\a(-n_1)...\a(-n_r){\bf 1} &\text {if $r$ is even;}\\
\a(-n_1)...\a(-n_{r-1}){\bf 1} & \text{if $r$ is odd;}\\
\end{cases}
\end{equation}
and
\begin{equation}
w=
\begin{cases}
e^{m\a}+e^{-m\a} &\text {if $r$ is even;}\\
\a(-n_r)(e^{m\a}-e^{-m\a}) &\text {if $r$ is odd.}\\
\end{cases}
\end{equation}
By Lemma \ref{sl4.5}, we have $v_{-1}w=u +u'$ where $u'\in
V_L^+(m)$ and $l_{\a}(u')<l_{\a}(u)=r$. Since $E^m, \a(-n_r)F^m\in
M(1)^++C_2(V_L^+)$, it implies that $u +u'\in  M(1)^+ +
C_2(V_L^+).$ By the induction hypothesis, we can conclude that $u
\in M(1)^+ + C_2(V_L^+).$ \qed

\smallskip

\begin{lem}\label{l4.13}
For any odd positive integer $m$, $V_L^+(m)\subset
V_L^+(1)+C_2(V_L^+).$
\end{lem}

\smallskip

\begin{thm}
$V_L^+/C_2(V_L^+)$ is spanned by $M(1)^+ +V_L^+(1)+C_2(V_L^+)$.
\end{thm}

\pf This follows from Lemmas \ref{l4.12}, \ref{l4.13} and the fact
that $$V_L^+=M(1)^+\oplus \oplus_{m=1}^{\infty}V_L^+(m). \ \
\qed$$
\medskip
Recall that $J={1\over 4k^2}\a(-1)^4{\bf 1}-{1\over
k}\a(-3)\a(-1){\bf 1}+{3\over 4k}\a(-2)^2{\bf 1}$.
\begin{prop} {[DN1]}
The vertex operator algebra $M(1)^+$ is spanned by
$$L(-m_1)...L(-m_s)J_{-n_1}...J_{-n_t}{\bf 1}$$
where $m_1\geq m_2\geq ...\geq m_s\geq 2$, $n_1\geq n_2\geq
...\geq n_t \geq 1$.
\end{prop}

\smallskip

\begin{Notation} Let $s,t\in \Z_{\geq 0}$, $v\in V_L^+$. We define
\begin{eqnarray*}
L(-2)^s&:=&(L(-2))^s\\
J_{-1}^t&:=& (J_{-1})^t.\\
\end{eqnarray*}
\end{Notation}

\smallskip

\begin{thm}\label{s413}
$(M(1)^++C_2(V_L^+))/C_2(V_L^+)$ is spanned by
$L(-2)^sJ_{-1}^t{\bf 1}+C_2(V_L^+)$ where $s, t\geq 0.$
\end{thm}

\smallskip

\begin{prop} {[DN2]} \ \

1) If $k$ is not  a perfect square, $V_L^+(1)$ is spanned by the
vectors
$$L(-n_1)L(-n_2)...L(-n_r)E, \ \ n_1\geq n_2\geq n_2\geq ...\geq n_t\geq 1.$$

\ \

2) If $k$ is a perfect square, $V_L^+(1)$ is spanned by the
vectors
$$L(-m_1)...L(-m_s)J_{-n_1}...J_{-n_t}E,\ \ m_1\geq m_2\geq ...\geq m_s\geq 1,\ \ n_1\geq n_2\geq ...\geq n_t \geq 1.$$
\end{prop}

\smallskip

\begin{thm}\label{st415} \ \

1) If $k$ is not a perfect square,
$(V_L^+(1)+C_2(V_L^+))/C_2(V_L^+)$ is spanned by
$L(-2)^sE+C_2(V_L^+)$ where $s\geq 0$.

\ \

2) If $k$ is a perfect square, $(V_L^+(1)+C_2(V_L^+))/C_2(V_L^+)$
is spanned by $L(-2)^sJ_{-1}^tE+C_2(V_L^+)$ where $s,t\geq 0$.

\end{thm}

\bigskip

\section{\bf{ The cofiniteness $C_2$ condition of $V_L^+$}}

In this section, we show that {\em $V_L^+$ satisfies the $C_2$
condition}. In particular, we prove that
$\left(V_L^+(1)+C_2(V_L^+)\right)/C_2(V_L^+)$ and
$\left(M(1)^++C_2(V_L^+)\right)/C_2(V_L^+)$ are finite
dimensional. We begin with showing that there is $\b \in \C$ such that for any $n\in \Z_{>0}$, $J_{-1}^nE\equiv \b ^n L(-2)^{2n}E \ \
\mod \ \ C_2(V_L^+).$ Consequently, $V_L^+(1)+C_2(V_L^+)$ is
spanned by $L(-2)^s E +C_2(V_L^+)$, $s\geq 0$. By using
information about the basis of $V_L^+(1)$, we are be able to prove
that $\left(V_L^+(1)+C_2(V_L^+)\right)/C_2(V_L^+)$ has finite
dimension. By taking the same approach, we can show that
$\left(M(1)^++C_2(V_L^+)\right)/C_2(V_L^+)$ has finite dimension.
Indeed, we have the following theorem.

\smallskip

\begin{thm}\ \

1)\begin{eqnarray*} \frac{V_L^+(1)+C_2(V_L^+)}{C_2(V_L^+)}&=&\<L(-2)^iE+C_2(V_L^+)\ \ |\ \ i=0,1,2\>\ \ \mbox{and}\\
\frac{M(1)^++C_2(V_L^+)}{C_2(V_L^+)}&=&\<L(-2)^m{\bf
1}+C_2(V_L^+),
J+C_2(V_L^+), J_{-1}^2{\bf 1}+C_2(V_L^+),\\
& & \ \ \ L(-2)J+C_2(V_L^+)\ \ |\ \ 0\leq m\leq 2k+2 \>.
\end{eqnarray*}

\ \

2) $V_L^+$ satisfies the $C_2$ condition.
\end{thm}

\smallskip

 For $m\in \Z_{>0}$, we set
$$V_L^+(m)=\oplus_{n\geq m^2k}V_L^+(m,n)$$ where $V_L^+(m,n)$ is the weight $n$ subspace of $V_L^+(m)$.
It is easy to see that the following elements form the bases of
$V_L^+(1,k+4)$ and $V_L^+(1,k+3)$, respectively.

\medskip

\underline{Basis of $V_L^+(1, k+4)$}
\begin{eqnarray*}
a_1&=& \a(-1)^4 E\\
a_2&=& \a(-2)\a(-1)^2F\\
a_3&=& \a(-3)\a(-1)E\\
a_4&=& \a(-4)F\\
a_5&=& \a(-2)^2E\\
\end{eqnarray*}

\underline{Basis of $V_L^+(1, k+3)$}
\begin{eqnarray*}
b_1&=& \a(-1)^3F\\
b_2&=& \a(-2)\a(-1)E\\
b_3&=& \a(-3)F.\\
\end{eqnarray*}

\smallskip

\begin{lem}\label{l5.1}
The vectors $L(-4)E$, $L(-1)b_i$ $(i=1,2,3)$, and $L(-2)^2E$ form
a basis of $V_L^+(1, k+4)$.
\end{lem}
\pf The table 1 in the appendix gives explicit expressions of
$L(-4)E$, $L(-1)b_i$ $(i=1,2,3)$ and $L(-2)^2E$ in terms of $a_j$,
$j=1,...,5$. If we denote the table 1 by $5\times 5$ matrix $A$,
then $\det (A)=-{9-40k+16k^2\over 16k^2}$. Since $9-40k+16k^2=0$
only when $k={1\over 4}, {9\over 4}$, it implies that $A$ is a
nonsingular matrix. Therefore, we conclude that $L(-4)E$,
$L(-1)b_i$ $(i=1,2,3)$, and $L(-2)^2E$ form a basis for
$V_L^+(1,k+4)$.

\smallskip

\begin{cor}\label{c5.2} $J_{-1}E\equiv \b L(-2)^2E \ \ \mod \ \ C_2(V_L^+)$ where
$\b= \frac{64k^2-16k-18}{(4k-1)(4k-9)}.$
\end{cor}
\pf Set
\begin{eqnarray*}
A_1&=& L(-2)^2E\\
A_2&=& L(-4)E\\
A_3&=& L(-1)\a(-1)^3F\\
A_4&=& L(-1)\a(-2)\a(-1)E\\
A_5&=& L(-1)\a(-3)F.\\
\end{eqnarray*}
The table 2 in the appendix represents the inverse matrix of the
matrix $A$. Since $$J_{-1}E={1\over 4k^2}a_1+{6\over
k}a_2+(12-{1\over k})a_3+(8k-17)a_4+(6+{3\over 4k})a_5,$$ we can
also think of $J_{-1}E$ as a vector $u=[{1\over 4k^2},{6\over
k},12-{1\over k},8k-17,6+{3\over 4k}]$. By multiplying the vector
$u$ and the matrix $A^{-1}$ together, we obtain that
$J_{-1}E\equiv \b L(-2)^2E \ \ \mod \ \ C_2(V_L^+)$. \qed

\smallskip

\begin{lem}\label{l5.3} For any positive integers $m,n$,

\ \ 1) $J_{-1}^nE \equiv \b ^nL(-2)^{2n}\ \ \mod\ \
C_2(V_L^+)$.

\ \ 2) $L(-2)^mJ_{-1}^nE\equiv \b ^nL(-2)^{2n+m}E\ \
\mod\ \ C_2(V_L^+)$.
\end{lem}

\pf Part 1) we will prove by induction on $n$. For $n=1$, it
follows from Corollary \ref{c5.2}. Now, suppose that
$J_{-1}^nE\equiv \left( \b\right)^nL(-2)^{2n}E\ \ \mod\ \
C_2(V_L^+)$. By using commutativity and associativity of
$V_L^+/C_2(V_L^+)$, we obtain the following:
\begin{eqnarray*}
J_{-1}^{n+1}E& \equiv & J_{-1}\left(\b\right)^nL(-2)^{2n}E\ \ \mod\ \ C_2(V_L^+) \\
&=&\left(\b\right)^nJ_{-1}\underbrace{\omega_{-1}...\omega_{-1}}_{2n\ \ \mbox{terms}}E\ \ \mod\ \ C_2(V_L^+)\\
&\equiv& \left(\b\right)^n \omega_{-1}...\omega_{-1}J_{-1}E \ \mod\ \ C_2(V_L^+)\\
&\equiv& \left(\b\right)^n \omega_{-1}...\omega_{-1}\b L(-2)^2E \ \mod\ \ C_2(V_L^+)\ \ (\text{by Corollary \ref{c5.2}})\\
&=& \left(\b\right)^{n+1} L(-2)^{2(n+1)}E \ \mod\ \ C_2(V_L^+).\\
\end{eqnarray*}
Part 2) it follows immediately from Part 1). \qed

\smallskip

\begin{cor} $(V_L^+(1)+C_2(V_L^+))/C_2(V_L^+)=\< \ \ L(-2)^sE+C_2(V_L^+)\ \ |\ \ s\geq 0 \ \ \>$.
\end{cor}
\pf This follows from Theorem \ref{st415} and Lemma
\ref{l5.3}.\qed

\medskip

The following elements are bases of $V_L^+(1, k+5)$, and
$V_L^+(1,k+3)$:

\smallskip

\underline{ Basis of $V_L^+(1, k+5)$}
\begin{eqnarray*}
B_1&=& \a(-5)F\\
B_2&=& \a(-4)\a(-1)E\\
B_3&=& \a(-3)\a(-2)E\\
B_4&=& \a(-3)\a(-1)^2F\\
B_5&=& \a(-2)^2\a(-1)F\\
B_6&=& \a(-2)\a(-1)^3E\\
B_7&=& \a(-1)^5F^2\\
\end{eqnarray*}

\medskip

\underline{ Basis of $V_L^+(1, k+3)$}
\begin{eqnarray*}
C_1&=&\a(-3)F\\
C_2&=&\a(-2)\a(-1)E\\
C_3&=&\a(-1)^3F.\\
\end{eqnarray*}

\smallskip

\begin{thm}{[DN2]} The vectors $L(-1)B_i$ $(i=1,...,7)$, $L(-3)C_j$ $(j=1,2,3)$ and $(\a(-1)^4{\bf 1})_{-3}E$ form a basis of $V_L^+(1, k+6)$.
\end{thm}

\smallskip

\begin{cor}\ \

1) $V_L^+(1,k+6)$ is a subset of $C_2(V_L^+)$.

2) $L(-2)^3E\in C_2(V_L^+)$.
\end{cor}

\smallskip

\begin{thm} $(V_L^+(1)+C_2(V_L^+))/C_2(V_L^+)=\< \ \ L(-2)^iE+C_2(V_L^+)\ \ |\ \ i=0,1,2\ \ \>. $
\end{thm}
\pf This follows from the fact that $L(-2)^3E\in C_2(V_L^+)$ and
$C_2(V_L^+)$ is an ideal of $V_L^+$ under $(-1)^{st}$ product.
\qed

\smallskip

Next, we show that $(M(1)^++C_2(V_L^+))/C_2(V_L^+)$ has finite
dimension. We set $$M(1)^+=\bigoplus_{i=0}^{\infty}M(1,i)^+$$
where $M(1,i)^+$ is the weight $i$ subspace of $M(1)^+$. It is
easy to see that $M(1,8)^+$, $M(1,7)^+$ and $M(1,5)^+$ have the
following basis elements:

\medskip

\underline{ Basis of $M(1,8)^+$}
\begin{eqnarray*}
c_1&=& \a(-1)^8{\bf 1}\\
c_2&=& \a(-3)\a(-1)^5{\bf 1}\\
c_3&=& \a(-5)\a(-1)^3{\bf 1}\\
c_4&=& \a(-7)\a(-1){\bf 1}\\
c_5&=& \a(-6)\a(-2){\bf 1}\\
c_6&=& \a(-5)\a(-3){\bf 1}\\
c_7&=& \a(-4)\a(-2)\a(-1)^2{\bf 1}\\
c_8&=& \a(-4)^{2}{\bf 1}\\
c_9&=& \a(-3)^2\a(-1)^2{\bf 1}\\
c_{10}&=& \a(-3)\a(-2)^2\a(-1){\bf 1}\\
c_{11}&=& \a(-2)^2\a(-1)^4{\bf 1}\\
c_{12}&=& \a(-2)^4{\bf 1}.\\
\end{eqnarray*}

\medskip

\underline{ Basis of $M(1,7)^+$}
\begin{eqnarray*}
\alpha_1&=&\a(-2)\a(-1)^5{\bf 1}\\
\alpha_2&=&\a(-4)\a(-1)^3{\bf 1}\\
\alpha_3&=&\a(-6)\a(-1){\bf 1}\\
\alpha_4&=&\a(-5)\a(-2){\bf 1}\\
\alpha_5&=&\a(-4)\a(-3){\bf 1}\\
\alpha_6&=&\a(-3)\a(-2)\a(-1)^2{\bf 1}\\
\alpha_7&=&\a(-2)^3\a(-1){\bf 1}.\\
\end{eqnarray*}

\medskip

\underline{ Basis of $M(1,5)^+$}
\begin{eqnarray*}
\beta_1&=&\a(-2)\a(-1)^3{\bf 1}\\
\beta_2&=&\a(-4)\a(-1){\bf 1}\\
\beta_3&=&\a(-3)\a(-2)^3{\bf 1}.\\
\end{eqnarray*}

\smallskip

\begin{lem}\label{l5.8} $M(1,8)^+$ is spanned by $L(-1)\alpha_i$ $(i=1,...,7)$, $L(-3)\beta_j$ $(j=1,2,3)$, $L(-2)^2J$, and $L(-2)^4{\bf 1}$.
\end{lem}

\pf For $i\in\{1,...,7\}$, we set $C_i=L(-1)\alpha_i$. We also set
$C_8=L(-3)\beta_1$, $C_9=L(-3)\beta_2$, $C_{10}=L(-2)^2J$,
$C_{11}=L(-2)^4{\bf 1}$, and $C_{12}=L(-3)\beta_3$. The table 3 in
the appendix gives explicit expressions of $C_i$ $(i=1,...,12)$ in
term of $c_j$ $(j=1,...,12)$. If we denote this table by $12\times
12$ matrix $B$, then the determinant of the matrix $B$ is
$-\frac{36315}{128k^8}$. Thus, $B$ is a nonsingular matrix and
$C_i$ $(i=1,...,12)$ span $M(1,8)^+$.

\smallskip

\begin{cor}\label{c5.8} $J_{-1}J\equiv \rho L(-2)^2J +\sigma L(-2)^4{\bf 1}\ \ \mod\ \ C_2(V_L^+)$ where $\rho=3.28+0.098k^{-1}$, and $\sigma=2.87-0.39k^{-1}$.
\end{cor}
\pf Since
$J_{-1}J=\frac{1}{16k^4}c_1+\frac{11}{2k^3}c_2+\frac{12}{k^2}c_3+\frac{558}{k}c_4-\frac{87}{k}c_5+\frac{186}{k}c_6+\frac{90}{k^2}c_7+\frac{72}{k}c_8+\frac{43}{k^2}c_9+\frac{117}{2k^2}c_{10}+\frac{51}{8k^3}c_{11}+\frac{105}{16k^2}c_{12}$,
we can also think of $J_{-1}J$ as a vector
$$w=[\frac{1}{16k^4},\frac{11}{2k^3},\frac{12}{k^2},\frac{558}{k},-\frac{87}{k},\frac{186}{k},\frac{90}{k^2},\frac{72}{k},\frac{43}{k^2},\frac{117}{2k^2},\frac{51}{8k^3},\frac{105}{16k^2}
].$$ The table 4 and table 5 in the appendix represents an inverse matrix of
$B$. By multiplying the vector $w$ and the matrix $B^{-1}$
together, we then obtain that $J_{-1}J\equiv \rho L(-2)^2J +\sigma
L(-2)^4{\bf 1}\ \ \mod\ \ C_2(V_L^+)$ where
$\rho=3.06+0.098k^{-1}$, and $\sigma=3.73-0.39k^{-1}$. \qed

\medskip

\begin{lem}{[DN1]} The vectors $L(-1)M(1,9)^+$, $L(-3)M(1,7)^+$,
$(\a(-1)^4{\bf 1})_{-3}\a(-1)^4{\bf 1}$, and $L(-2)^5{\bf 1}$ span
$M(1,10)^+$.
\end{lem}

\smallskip

\begin{rem}{[DN1]} Only $L(-2)^5{\bf 1}$ involves the vector
$\a(-1)^{10}{\bf 1}$.
\end{rem}
\smallskip

\begin{cor} $(M(1,10)^++C_2(V_L^+))/C_2(V_L^+)=\<\ \ L(-2)^5{\bf
1}+C_2(V_L^+)\>.$ \end{cor}

\smallskip

\begin{cor}\label{c514}\ \

1) $L(-2)^5{\bf 1}\equiv \frac{1}{(4k)^5}\a(-1)^{10}{\bf 1}\ \
\mod \ \ C_2(V_L^+)$.

2) $L(-2)J_{-1}J\equiv 16 L(-2)^5{\bf 1} \ \ \mod \ \ C_2(V_L^+).$

3) $L(-2)^3 J \equiv 4L(-2)^5{\bf 1}\ \ \mod \ \ C_2(V_L^+).$
\end{cor}

\pf For 1), it follows from the fact that
\begin{eqnarray*}
L(-2)^4{\bf
1}&=&\frac{1}{(4k)^4}c_1+\frac{6}{k^3}c_2+\frac{3}{2k^2}c_3+\frac{15}{2k}c_4+\frac{9}{2k}c_6+\frac{5}{4k^2}c_9\
\ \mbox{and}\\
L(-2)\a(-1)^8{\bf 1}&=&8\a(-3)\a(-1)^7{\bf
1}+\frac{1}{4k}\a(-1)^{10}{\bf 1}.\end{eqnarray*}

For 2), we first recall that
\begin{eqnarray*}
J_{-1}J&=&\frac{1}{16k^4}c_1+\frac{11}{2k^3}c_2+\frac{12}{k^2}c_3+\frac{558}{k}c_4-\frac{87}{k}c_5+\frac{186}{k}c_6\\
& &+\frac{90}{k^2}c_7+\frac{72}{k}c_8+\frac{43}{k^2}c_9+\frac{117}{2k^2}c_{10}+\frac{51}{8k^3}c_{11}+\frac{105}{16k^2}c_{12}.\\
\end{eqnarray*}
Since $L(-2)\a (-1)^8{\bf 1}\equiv \frac{1}{4k}\a(-1)^{10}{\bf 1}
\ \ \mod \ \ C_2(V_L^+)$, this implies that $L(-2)J_{-1}J\equiv 16
L(-2)^5{\bf 1} \ \ \mod \ \ C_2(V_L^+)$.

For 3), it follows from the fact that
\begin{eqnarray*}
L(-2)^2J&=&\frac{1}{64k^4}c_1+\frac{9}{16k^3}c_2+\frac{3}{2k^2}c_3-\frac{15}{k}c_4+\frac{12}{k}c_5-\frac{9}{k}c_6+\frac{3}{2k^2}c_7+\frac{6}{k}c_8\\
& &+\frac{2}{k^2}c_9+\frac{3}{8k^2}c_{10}+\frac{3}{64k^3}c_{11},\\
L(-2)\a(-1)^8{\bf 1}&\equiv&\frac{1}{4k}\a(-1)^{10}{\bf 1}\ \ \mod
\ \ C_2(V_L^+),\ \ \mbox{and}\ \ 1). \ \ \qed\end{eqnarray*}

\smallskip

\begin{lem}\label{l515} $J_{-1}^2J\equiv \g L(-2)^6{\bf 1}\ \ \mod \ \
C_2(V_L^+)$ where $\g=16\rho+4\sigma$.
\end{lem}

\pf Recall from Corollary \ref{c5.8} that $J_{-1}J\equiv \rho
L(-2)^2J +\sigma L(-2)^4{\bf 1}\ \ \mod\ \ C_2(V_L^+)$. By
commutative, associative laws of $V_L^+/C_2(V_L^+)$, and Corollary
\ref{c514}, we have
\begin{eqnarray*}
J_{-1}J_{-1}J&\equiv& \rho J_{-1}L(-2)^2J +\sigma J_{-1}L(-2)^4{\bf 1}\ \ \mod\ \ C_2(V_L^+)\\
&\equiv& \rho L(-2)^2J_{-1}J+\sigma L(-2)^4 J\ \ \mod \ \
C_2(V_L^+)\\
&\equiv& 16\rho L(-2)^6{\bf 1}+4\sigma L(-2)^6{\bf 1} \ \ \mod \ \
C_2(V_L^+)\\
&=& \g L(-2)^6{\bf 1} \ \ \mod \ \ C_2(V_L^+). \qed
\end{eqnarray*}

\begin{lem}\label{l516} For $n\geq 4$, $J_{-1}^n{\bf 1}\equiv \g
4^{n-3}L(-2)^{2n}{\bf 1} \ \ \mod \ \ C_2(V_L^+).$ \end{lem}

\pf For $n=4$, this follows from Lemma \ref{l515} and Corollary
\ref{c514} 3) and commutative law of $V_L^+/C_2(V_L^+)$.

Suppose $J_{-1}^n{\bf 1}\equiv \g 4^{n-3}L(-2)^{2n}{\bf 1}\ \ \mod
\ \ C_2(V_L^+)$. By Corollary \ref{c514} 3) and commutativity of
$V_L^+/C_2(V_L^+)$, we have
\begin{eqnarray*}
J_{-1}^{n+1}{\bf 1}&\equiv& \g 4^{n-3}J_{-1}L(-2)^{2n}{\bf 1}\ \
\mod \ \
C_2(V_L^+)\\
&\equiv& \g 4^{n-3}L(-2)^{2n-3}L(-2)^3 J\ \ \mod \ \ C_2(V_L^+)\\
&\equiv& \g 4^{n-3}4L(-2)^{2n-3}L(-2)^5{\bf 1}\ \ \mod \ \
C_2(V_L^+)\\
&\equiv& \g 4^{(n+1)-3}L(-2)^{2(n+1)}{\bf 1}\ \ \mod \ \
C_2(V_L^+). \\
\end{eqnarray*}
Therefore, $J_{-1}^n{\bf 1}\equiv \g 4^{n-3}L(-2)^{2n}{\bf 1} \ \
\mod \ \ C_2(V_L^+)$ for all $n\geq 4$.\qed

\smallskip

\begin{cor}\label{l5.9} For $n\geq 3$, $J_{-1}^n{\bf 1} \equiv \g 4^{n-3} L(-2)^{2n}{\bf
1}\ \ \mod \ \ C_2(V_L^+)$.
\end{cor}

\smallskip

\begin{cor}\label{c5.10} For $m\in \Z_+$, $n\geq 3$, $L(-2)^m J_{-1}^n{\bf 1}\equiv \g
4^{n-3} L(-2)^{m+2n}{\bf 1} \ \ \mod \ \ C_2(V_L^+)$.

\end{cor}

\smallskip

\begin{cor} $(M(1)^++C_2(V_L^+))/C_2(V_L^+)$ is spanned by $L(-2)^s{\bf 1}+C_2(V_L^+)$, $J+C_2(V_L^+),$ $J_{-1}J+C_2(V_L^+)$, $L(-2)J+C_2(V_L^+)$ where $s\geq 0$.
\end{cor}
\pf This follows from Theorem \ref{s413}, Corollary \ref{c5.8},
Corollary \ref{l5.9}, and Corollary \ref{c5.10}. \qed

\smallskip

The following elements are bases of $V_L^+(2,4k+6)$,
$V_L^+(2,4k+5)$ and $V_L^+(2,4k+3)$:

\smallskip

\underline{ Basis of $V_L^+(2,4k+6)$}

\begin{eqnarray*}
g_1&=&\a(-6)F^2\\
g_2&=&\a(-5)\a(-1)E^2\\
g_3&=&\a(-4)\a(-2)E^2\\
g_4&=&\a(-4)\a(-1)^2F^2\\
g_5&=&\a(-3)^2E^2 \\
g_6&=&\a(-3)\a(-2)\a(-1)F^2\\
g_7&=&\a(-3)\a(-1)^3E^2\\
g_8&=&\a(-2)^3F^2\\
g_9&=&\a(-2)^2\a(-1)^2E^2 \\
g_{10}&=&\a(-2)\a(-1)^4F^2\\
g_{11}&=&\a(-1)^6E^2
\end{eqnarray*}

\underline{ Basis of $V_L^+(2, 4k+5)$}
\begin{eqnarray*}
f_1&=& \a(-5)F^2\\
f_2&=& \a(-4)\a(-1)E^2\\
f_3&=& \a(-3)\a(-2)E^2\\
f_4&=& \a(-3)\a(-1)^2F^2\\
f_5&=& \a(-2)^2\a(-1)F^2\\
f_6&=& \a(-2)\a(-1)^3E^2\\
f_7&=& \a(-1)^5F^2\\
\end{eqnarray*}

\underline{ Basis of $V_L^+(2, 4k+3)$}
\begin{eqnarray*}
h_1&=&\a(-3)F^2\\
h_2&=&\a(-2)\a(-1)E^2\\
h_3&=&\a(-1)^3F^2.\\
\end{eqnarray*}

\smallskip

\begin{thm}\label{t5.13}
$V_L^+(2, 4k+6)$ is spanned by $L(-1)f_i$ $(i=1,...,7)$,
$L(-3)h_j$ $(j=1,2,3)$ and $\left(\a(-1)^4\bf 1\right)_{-3}E^2$.
\end{thm}
\pf Set
\begin{eqnarray*}
G_1&=& L(-1)\a(-5)F^2\\
G_2&=& L(-1)\a(-4)\a(-1)E^2\\
G_3&=& L(-1)\a(-3)\a(-2)E^2\\
G_4&=& L(-1)\a(-3)\a(-1)^2F^2\\
G_5&=& L(-1)\a(-2)^2\a(-1)F^2\\
G_6&=& L(-1)\a(-2)\a(-1)^3E^2\\
G_7&=& L(-1)\a(-1)^5F^2\\
G_8&=& L(-3)\a(-3)F^2\\
G_9&=& L(-3)\a(-2)\a(-1)E^2\\
G_{10}&=& L(-3)\a(-1)^3F^2\\
G_{11}&=& \left( \a(-1)^4{\bf 1}\right)_{-3}E^2.\\
\end{eqnarray*}
The table 6 in the appendix gives explicit expressions of $G_i$
$(i=1,...,11)$ in terms of $g_j, j=l,...,11$. If we denote this
table by $11\times 11$ matrix $C$, then $det(C)=-{24\over
k^2}({1536k^4-2592k^3+1072k^2-58k+15})$. Note that $C$ is
nonsingular if $k$ is a positive integer. Therefore, $V_L^+(2,
4k+6)$ is spanned by $L(-1)f_i$ $(i=1,...,7)$, $L(-3)h_j$
$(j=1,2,3)$ and $\left(\a(-1)^4\bf 1\right)_{-3}E^2$. \qed

\smallskip

\begin{cor}\label{c5.14} $V_L^+(2, 4k+6)$ is contained in $C_2(V_L^+)$.
\end{cor}

\smallskip

\begin{lem}
$(M(1)^++C_2(V_L^+))/C_2(V_L^+)$ is spanned by $ J+C_2(V_L^+),\ \
J_{-1}J+C_2(V_L^+),\ \ L(-2)J+C_2(V_L^+),\ \ L(-2)^m{\bf
1}+C_2(V_L^+)$ where $0\leq m\leq 2k+2.$
\end{lem}
\pf We will show that $L(-2)^{2k+3}{\bf 1}\in C_2(V_L^+)$. Suppose
the contrary. By Corollaries \ref{c514}, \ref{c5.10}, we conclude
that $L(-2)^mJ_{-1}^n{\bf 1}\not\in C_2(V_L^+)$ for all $m, n \in
\Z_{+}$ such that $m+2n=2k+3$. Moreover, we have
$$M(1,4k+3)^+=\<L(-2)^{2k+3}{\bf 1}\>\oplus \left(C_2(V_L^+)\cap
M(1,4k+3)^+\right).$$

Observe that $E_{-2k-7}E=p_6(\a)\o e^{2\a}+p_6(-\a)\o
e^{-2\a}+p_{4k+6}(\a)+p_{4k+6}(-\a)$. Since $p_6(\a)\o
e^{2\a}+p_6(-\a)\o e^{-2\a}\in C_2(V_L^+)$ ( by Corollary
\ref{c5.14}), this implies that $p_{4k+6}(\a)+p_{4k+6}(-\a)\in
C_2(V_L^+).$ Since $p_{4k+6}(\a)+p_{4k+6}(-\a)$ is also an element
of $M(1)^+$, we can write $p_{4k+6}(\a)+p_{4k+6}(-\a)$ in the
following way:
\begin{eqnarray*}
(*)\ \ p_{4k+6}(\a)+p_{4k+6}(-\a)&=&\sum\l_iL(-m_1)...L(-m_i){\bf 1}\\
& &\ \ +\sum \r_{j,l}L(-n_1)...L(-n_{j})J_{-p_1}...J_{-p_l}{\bf 1}\\
& &\ \ + \sum\b_t J_{-q_1}...J_{-q_t}{\bf 1}.\ \ \
\end{eqnarray*}
Here,

1) $m_1\geq m_2\geq ...\geq m_i\geq 2$ and there is $1\leq f\leq
i$ such that $m_f\geq 3$. Moreover, $m_1+...+m_i=4k+6$.

2) $n_1\geq n_2 \geq ... \geq n_j\geq 2$, $p_1\geq ...\geq p_l\geq
1$ and there is $1\leq g\leq j$ such that $n_g\geq 3$ or there is
$1\leq h\leq l$ such that $p_h\geq 2$. Furthermore,
$n_1+...+n_j+p_1+...+p_l+3l=4k+6$.

3) $q_1\geq ...\geq q_t\geq 1$ and there is $1\leq s \leq t$ such
that $q_s\geq 2$. Moreover, $q_1+...+q_t+3t=4k+6.$
\smallskip

We observe that in the spanning set of $M(1)^+_{4k+6}$ only
$L(-2)^{2k+3}{\bf 1}$ and $L(-2)^mJ_{-1}^n{\bf 1}$ involve with
the term $\a(-1)^{4k+6}{\bf 1}$. Here, $m,n\in \Z_{+}$ such that
$m+2n=2k+3$. Since $L(-2)^{2k+3}{\bf 1},L(-2)^mJ_{-1}^n{\bf 1}$ do
not occur on the right side of $(*)$, it implies that
$\a(-1)^{4k+6}{\bf 1}$ is a linear combination of vectors that
have length of $\a$ less than $4k+6$. This is impossible. So, we
have a contradiction. Thus, $L(-2)^{2k+3}{\bf 1}\in C_2(V_L^+)$.
Moreover, \begin{eqnarray*} & &(M(1)^++C_2(V_L^+))/C_2(V_L^+)\\
&=&\< J+C_2(V_L^+),J_{-1}J+C_2(V_L^+),L(-2)J+C_2(V_L^+),L(-2)^m{\bf 1}+C_2(V_L^+)|\\
& & \ \ \ 0\leq m\leq 2k+2\>. \ \ \qed\\
\end{eqnarray*}

\smallskip

\begin{thm} $V_L^+$ satisfies the $C_2$ condition.
\end{thm}

\smallskip

\begin{cor} $V_L^+$ satisfies the cofiniteness $C_n$ condition for
all $n\in \Z_{>0}$ such that $n\geq 2$.

\end{cor}

\smallskip

\begin{cor}
For $n\in \Z_{\geq 0}$, $A_n(V_L^+)$ is finite dimensional.
\end{cor}

\clearpage

\section{Appendix}
\begin{table}[ht]
\begin{center}
\caption{}
\begin{tabular}[ht]{|| c | c | c | c | c | c ||}
\hline
& & & & &\\
& $a_1$ &$a_2$ &$a_3$ &$a_4$ &$a_5$ \\
& & & & &\\
\hline
& & & & & \\
$L(-2)^2E$& $1/( 4k)^2$& $1/(2k)$&$1/(2k)$& 2& 1\\
& & & & & \\
\hline
& & & & & \\
$L(-4) E$& 0& 0& $1/(2k)$ &1 &$1/(4k)$\\
& & & & &\\
\hline
& & & & &\\
$L(-1)\a (-1)^3F$& 1& 3& 0& 0& 0\\
& & & & & \\
\hline
& & & & & \\
$L(-1)\a (-2)\a (-1)E$& 0&1& 2& 0&1\\
& & & & & \\
\hline
& & & & & \\
$L(-1)\a (-3)F$& 0& 0& 1& 3&0\\
& & & & & \\
\hline
\end{tabular}
\end{center}
\end{table}
\medskip
\begin{table}[ht]
\begin{center}
\caption{}
\begin{tabular}[ht]{||c|c|c|c|c|c||}
\hline
& & & & & \\
& $A_1$ & $A_2$ & $A_3$ & $A_4$ & $A_5$\\
& & & & & \\
\hline
& & & & & \\
$a_1$ & \small{${48k^2\over(4k-1)(4k-9)}$}& \small{${-24(16k-3)k\over (4k-1)(4k-9)}$} & \small{${2(8k^2-20k+3)\over (4k-1)(4k-9)}$} & \small{$-{6(8k^2-16k+3)\over (4k-1)(4k-9)}$} & \small{$24k\over 4k-9$}\\
& & & & & \\
\hline
& & & & & \\
$a_2$ &  \small{${-16k^2\over(4k-1)(4k-9)}$}& \small{${8(16k-3)k\over (4k-1)(4k-9)}$} & \small{${1\over (4k-1)(4k-9)}$} & \small{${2(8k^2-16k+3)\over (4k-1)(4k-9)}$} & \small{$-{8k\over 4k-9}$} \\
& & & & & \\
\hline
& & & & & \\
$a_3$& \small{${12k\over (4k-1)(4k-9)}$}& \small{$-{3(16k^2-8k+3)\over (4k-1)(4k-9)}$}& \small{$-{3\over 4k(4k-1)(4k-9)}$}& \small{$-{72k-27\over 12k(4k-1)(4k-9)}$}& \small{${4k-3\over (4k-9)}$}\\
& & & & & \\
\hline
& & & & & \\
$a_4$&  \small{${-4k\over (4k-1)(4k-9)}$}& \small{${16k^2-8k+3\over (4k-1)(4k-9)}$} & \small{${1\over 4k(4k-1)(4k-9)}$}& \small{${8k-3\over 4k(4k-1)(4k-9)}$}& \small{$-{2\over 4k-9}$} \\
& & & & & \\
\hline
& & & & & \\
$a_5$&  \small{${8(2k-3)k\over (4k-1)(4k-9)}$}& \small{ $-{2(16k^2+12k-9)\over (4k-1)(4k-9)}$}& \small{$-{2k-3\over 2k(4k-1)(4k-9)}$}& \small{$-{(2k-3)(8k-3)\over 2k(4k-1)(4k-9)}$}& \small{${6\over 4k-9}$}\\
& & & & & \\
\hline
\end{tabular}
\end{center}
\end{table}
\clearpage
\begin{table}[ht]
\begin{center}
\caption{}
\begin{tabular}[ht]{||c|c|c|c|c|c|c|c|c|c|c|c|c||}
\hline
&&&&&&&&&&&&\\
&
$c_1$&$c_2$&$c_3$&$c_4$&$c_5$&$c_6$&$c_7$&$c_8$&$c_9$&$c_{10}$&$c_{11}$&$c_{12}$\\
&&&&&&&&&&&&\\
\hline
&&&&&&&&&&&&\\
$C_1$&0&2&0&0&0&0&0&0&0&0&5&0\\
&&&&&&&&&&&&\\
\hline
&&&&&&&&&&&&\\
$C_2$&0&0&4&0&0&0&3&0&0&0&0&0\\
&&&&&&&&&&&&\\
\hline
&&&&&&&&&&&&\\
$C_3$&0&0&0&6&1&0&0&0&0&0&0&0\\
&&&&&&&&&&&&\\
\hline
&&&&&&&&&&&&\\
$C_4$&0&0&0&0&5&2&0&0&0&0&0&0\\
&&&&&&&&&&&&\\
\hline
&&&&&&&&&&&&\\
$C_5$&0&0&0&0&0&4&0&3&0&0&0&0\\
&&&&&&&&&&&&\\
\hline
&&&&&&&&&&&&\\
$C_6$&0&0&0&0&0&0&3&0&2&2&0&0\\
&&&&&&&&&&&&\\
\hline
&&&&&&&&&&&&\\
$C_7$&0&0&0&0&0&6&0&0&0&0&0&1\\
&&&&&&&&&&&&\\
\hline
&&&&&&&&&&&&\\
$C_8$&0&0&2&0&0&0&3&0&0&0&{$1/ 2k$}&0\\
&&&&&&&&&&&&\\
\hline
&&&&&&&&&&&&\\
$C_9$&0&0&0&4&0&0&$1/2k$&1&0&0&0&0\\
&&&&&&&&&&&&\\
\hline
&&&&&&&&&&&&\\
$C_{10}$&$1/64k^4$&$9/16k^3$&$3/2k^2$&$-15/k$&$12/k$&$-9/k$&$3/2k^2$&$6/k$&$2/k^2$&$3/8k^2$&$3/64k^3$&0\\
&&&&&&&&&&&&\\
\hline
&&&&&&&&&&&&\\
$C_{11}$&$1/(4k)^4$&$6/k^3$&$3/2k^2$&$15/2k$& 0&$9/2k$&0&0&$5/4k^2$&0&0&0\\
&&&&&&&&&&&&\\
\hline
&&&&&&&&&&&&\\
$C_{12}$&0&0&0&0&3&2&0&0&0&$1/2k$&0&0\\
&&&&&&&&&&&&\\
\hline
\end{tabular}
\end{center}
\end{table}
\clearpage
\begin{table}[ht]
\begin{center}
\caption{}
\begin{tabular}[ht]{||c|c|c|c|c|c|c||}
\hline
&&&&&&\\
&
$C_1$&$C_2$&$C_3$&$C_4$&$C_5$&$C_6$\\
&&&&&&\\
\hline
&&&&&&\\
$c_1$&$-\frac{13392k}{269}$&$-\frac{73464k^2}{269}$&$\frac{1164704k^3}{1345}$&$-\frac{613408k^3}{1345}$&$\frac{323568k^3}{1345}$&$-\frac{91552k^2}{1345}$\\
&&&&&&\\
\hline
&&&&&&\\
$c_2$&$\frac{19}{538}$&$\frac{35k}{269}$&$-\frac{1204k^2}{807}$&$\frac{188k^2}{807}$&$-\frac{418k^2}{807}$&$-\frac{16k}{269}$\\
&&&&&&\\
\hline
&&&&&&\\
$c_3$&$\frac{25}{538k}$&$\frac{131}{269}$&$\frac{602k}{4035}$&$-\frac{94k}{4035}$&$\frac{209k}{4035}$&$\frac{8}{1345}$\\
&&&&&&\\
\hline
&&&&&&\\
$c_4$&$-\frac{25}{12912k^2}$&$-\frac{85}{8608k}$&$\frac{19573}{96840}$&$-\frac{3941}{96840}$&$\frac{3617}{193680}$&$-\frac{1}{4035k}$\\
&&&&&&\\
\hline
&&&&&&\\
$c_5$&$\frac{25}{2152k^2}$&$\frac{255}{4304k}$&$-\frac{3433}{16140}$&$\frac{3941}{16140}$&$-\frac{3617}{32280}$&$\frac{2}{1345k}$\\
&&&&&&\\
\hline
&&&&&&\\
$c_6$&$-\frac{125}{4304k^2}$&$-\frac{1275}{8608k}$&$\frac{3433}{6456}$&$-\frac{713}{6456}$&$\frac{3617}{12912}$&$-\frac{1}{269k}$\\
&&&&&&\\
\hline
&&&&&&\\
$c_7$&$-\frac{50}{807k}$&$-\frac{85}{269}$&$-\frac{2408k}{12105}$&$\frac{376k}{12105}$&$-\frac{836k}{12105}$&$-\frac{32}{4035}$\\
&&&&&&\\
\hline
&&&&&&\\
$c_8$&$\frac{125}{3228k^2}$&$\frac{425}{2152k}$&$-\frac{3433}{4842}$&$\frac{713}{4842}$&$-\frac{389}{9684}$&$\frac{4}{807k}$\\
&&&&&&\\
\hline
&&&&&&\\
$c_9$&$\frac{25}{538k}$&$\frac{255}{1076}$&$\frac{4637k}{4035}$&$\frac{3941k}{4035}$&$\frac{4453k}{8070}$&$\frac{1361}{2690}$\\
&&&&&&\\
\hline
&&&&&&\\
$c_{10}$&$\frac{25}{538k}$&$\frac{255}{1076}$&$-\frac{3433k}{4035}$&$-\frac{4129k}{4035}$&$-\frac{3617k}{8070}$&$\frac{8}{1345}$\\
&&&&&&\\
\hline
&&&&&&\\
$c_{11}$&$\frac{50}{269}$&$-\frac{14k}{269}$&$\frac{2408k^2}{4035}$&$-\frac{376k^2}{4035}$&$\frac{836k^2}{4035}$&$\frac{32k}{1345}$\\
&&&&&&\\
\hline
&&&&&&\\
$c_{12}$&$\frac{375}{2152k^2}$&$\frac{3825}{4304k}$&$-\frac{3433}{1076}$&$\frac{713}{1076}$&$-\frac{3617}{2152}$&$\frac{6}{269k}$\\
&&&&&&\\
\hline

\end{tabular}
\end{center}
\end{table}
\clearpage
\begin{table}[ht]
\begin{center}
\caption{}
\begin{tabular}[ht]{||c|c|c|c|c|c|c||}
\hline
&&&&&&\\
&
$C_7$&$C_8$&$C_9$&$C_{10}$&$C_{11}$&$C_{12}$\\
&&&&&&\\
\hline
&&&&&&\\
$c_1$& 0&$\frac{661872k^2}{1345}$&$-\frac{1465296k^3}{1345}$&$\frac{82432k^4}{1345}$&$\frac{14592k^4}{1345}$&$\frac{304384k^3}{1345}$\\
&&&&&&\\
\hline
&&&&&&\\
$c_2$&0&$-\frac{94k}{269}$&$\frac{482k^2}{269}$&$-\frac{32k^3}{807}$&$\frac{128k^3}{807}$&$\frac{72k^2}{269}$\\
&&&&&&\\
\hline
&&&&&&\\
$c_3$&0&$-\frac{1251}{2690}$&$-\frac{241k}{1345}$&$\frac{16k^2}{4035}$&$-\frac{64k^2}{4035}$&$-\frac{36k}{1345}$\\
&&&&&&\\
\hline
&&&&&&\\
$c_4$&0&$\frac{417}{21520k}$&$-\frac{3553}{64560}$&$-\frac{2k}{12105}$&$\frac{8k}{12105}$&$\frac{3}{2690}$\\
&&&&&&\\
\hline
&&&&&&\\
$c_5$&0&$-\frac{1251}{10760k}$&$\frac{3553}{10760}$&$\frac{4k}{4035}$&$-\frac{16k}{4035}$&$-\frac{9}{1345}$\\
&&&&&&\\
\hline
&&&&&&\\
$c_6$&0&$\frac{1251}{4304k}$&$-\frac{3553}{4304}$&$-\frac{2k}{807}$&$\frac{8k}{807}$&$\frac{9}{538}$\\
&&&&&&\\
\hline
&&&&&&\\
$c_7$&0&$\frac{834}{1345}$&$\frac{964k}{4035}$&$-\frac{64k^2}{12105}$&$\frac{256k^2}{12105}$&$\frac{48k}{1345}$\\
&&&&&&\\
\hline
&&&&&&\\
$c_8$&0&$-\frac{417}{1076k}$&$\frac{3553}{3228}$&$\frac{8k}{2421}$&$-\frac{32k}{2421}$&$-\frac{6}{269}$\\
&&&&&&\\
\hline
&&&&&&\\
$c_9$&0&$-\frac{1251}{2690}$&$-\frac{4517k}{2690}$&$\frac{16k^2}{4035}$&$-\frac{64k^2}{4035}$&$-\frac{2726k}{1345}$\\
&&&&&&\\
\hline
&&&&&&\\
$c_{10}$&0&$-\frac{1251}{2690}$&$\frac{3553k}{2690}$&$\frac{16k^2}{4035}$&$-\frac{64k^2}{4035}$&$\frac{2654k}{1345}$\\
&&&&&&\\
\hline
&&&&&&\\
$c_{11}$&0&$\frac{188k}{1345}$&$-\frac{964k^2}{1345}$&$\frac{64k^3}{4035}$&$-\frac{256k^3}{4035}$&$-\frac{144k^2}{1345}$\\
&&&&&&\\
\hline
&&&&&&\\
$c_{12}$&1&$-\frac{3753}{2152k}$&$\frac{10659}{2152}$&$\frac{4k}{269}$&$-\frac{16k}{269}$&$-\frac{27}{269}$\\
&&&&&&\\
\hline
\end{tabular}
\end{center}
\end{table}
\clearpage
\begin{table}[ht]
\begin{center}
\caption{}
\begin{tabular}[ht]{|| c|c|c|c|c|c|c|c|c|c|c|c||}
\hline
&&&&&&&&&&&\\
& $g_1$&$g_2$&$g_3$&$g_4$&$g_5$&$g_6$&$g_7$&$g_8$&$g_9$&$g_{10}$&$g_{11}$\\
&&&&&&&&&&&\\
\hline
&&&&&&&&&&&\\
$G_1$& 5& 2&0&0&0&0&0&0&0&0&0\\
&&&&&&&&&&&\\
\hline
&&&&&&&&&&&\\
$G_2$&0&4&1&2&0&0&0&0&0&0&0\\
&&&&&&&&&&&\\
\hline
&&&&&&&&&&&\\
$G_3$&0&0&3&0&2&2&0&0&0&0&0\\
&&&&&&&&&&&\\
\hline
&&&&&&&&&&&\\
$G_4$&0&0&0&3&0&2&2&0&0&0&0\\
&&&&&&&&&&&\\
\hline
&&&&&&&&&&&\\
$G_5$&0&0&0&0&0&4&0&1&2&0&0\\
&&&&&&&&&&&\\
\hline
&&&&&&&&&&&\\
$G_6$&0&0&0&0&0&0&2&0&3&2&0\\
&&&&&&&&&&&\\
\hline
&&&&&&&&&&&\\
$G_7$&0&0&0&0&0&0&0&0&0&5&2\\
&&&&&&&&&&&\\
\hline
&&&&&&&&&&&\\
$G_8$&3&0&0&0&2&$1/ 2k$&0&0&0&0&0\\
&&&&&&&&&&&\\
\hline
&&&&&&&&&&&\\
$G_9$&0&2&1&0&0&1&0&0&$1/2k$&0&0\\
&&&&&&&&&&&\\
\hline
&&&&&&&&&&&\\
$G_{10}$&0&0&0&3&0&0&2&0&0&$1/2k$&0\\
&&&&&&&&&&&\\
\hline
&&&&&&&&&&&\\
$G_{11}$&\small{$256k^3$}&\small{$192k^2$}&\small{$192k^2$}&\small{$48k$}&\small{$96k^2$}&\small{$96k$}&\small{$4$}&\small{$16k$}&\small{$6$}&0&0\\
&&&&&&&&&&&\\
\hline
\end{tabular}
\end{center}
\end{table}

\end{document}